\documentclass[10pt,reqno]{amsart}
\usepackage{latexsym}
\usepackage{amsmath}
\usepackage{setspace}
\usepackage{amsmath}

\newcommand{\lb}{\left(}

\newcommand{\vp}{\varphi}
\newcommand{\ve}{\varepsilon}

\newcommand{\rb}{\right)}
\newcommand{\PD}{\partial}

\newcommand{\Hc}{\mathcal{H}}

\newcommand{\Wc}{\mathcal{W}}

\newcommand{\R}{\rangle}

\newcommand{\La}{\langle}

\newcommand{\g}{\gamma}

\newcommand{\n}{\nabla}

\newcommand{\A}{\alpha}

\usepackage{amssymb,amsmath,amsfonts,amsthm}
\newtheorem{thm}{Theorem}
\newtheorem{lemma}{Lemma}

\newtheorem{defn}{Definition}

\theoremstyle{definition}

\begin{document}
\title[Inversion formulas]{A generalization of inversion formulas of Pestov and Uhlmann}
\author[V. P. Krishnan]{Venkateswaran P. Krishnan}
\address{Department of Mathematics,
Tufts University, 
Medford, MA 02155 USA.}
\email{venky.krishnan@tufts.edu}
\begin{abstract}
In this note, we give a generalization of the inversion formulas of Pestov-Uhlmann for the geodesic ray transform of functions and vector fields on simple $2$-dimensional manifolds of constant curvature. The inversion formulas given here hold for $2$-dimensional simple manifolds whose curvatures close to a constant.
\end{abstract} 
\maketitle
\section{Introduction}
Let $(M,\PD M,g)$ be a $C^{\infty}$ compact Riemannian manifolds with boundary. A variant of the classical Radon transform on Euclidean space is the {\it geodesic ray transform} on Riemannian manifolds defined as follows:
\[ I_{m}f(\g)=\int\limits_{0}^{l(\g)}\La f(\g(t)),\dot{\g}^{m}(t)\R dt,\]
where $\g:[0,l(\g)]\to M$ is a maximal geodesic parameterized by arc length and $m$ indicates the rank of the symmetric tensor field $f\in L^{2}(M)$.

We will be interested only in the cases $m=0$ (functions) and $m=1$ (vector fields) in this note and denote their geodesic ray transforms by $I_{0}$ and $I_{1}$ respectively.

The geodesic ray transform is not injective in general. One needs additional restrictions on the metric and one such restriction is to assume that the Riemannian manifold $(M,\PD M,g)$ is {\it simple} \cite{Sh1} defined as follows:
\begin{defn}
A compact Riemannian manifold with boundary is simple if 
\begin{enumerate}
\item[(a)] The boundary $\PD M$ is strictly convex: $\langle \n_{\xi}\nu,\xi\R < 0$ for $\xi \in T_{x}(\PD M)$ where $\nu$ is the unit inward normal to the boundary.
\item[(b)] The map $\exp_{x}:\exp_{x}^{-1}M\to M$ is a diffeomorphism for each $x\in M$.
\end{enumerate}
\end{defn}

It is known that on a simple Riemannian manifold, $I_{0}f$ uniquely determines $f$ and $I_{1}f$ uniquely determines the solenoidal component of $f$. For references to these works, we refer the book of Sharafutdinov \cite{Sh1}. Then, similar to the classical Radon inversion formula, it is natural to ask whether there exists explicit inversion formulas for a function or a vector field in terms of its geodesic ray transform. In general this is a hard problem and such formulas are known in only in special cases \cite{SH}.

In \cite{PU2}, Pestov and Uhlmann found Fredholm-type inversion formulas for the geodesic ray transform of functions and vector fields for simple 2-dimensional manifolds. These formulas become exact inversion formulas for $2$-dimensional manifolds of constant curvature, even when conjugate points are present along geodesics. 

A brief remark regarding notation. We use notation that is standard in integral geometry literature. Ours is consistent for the most part with \cite{Sh1}. In this note,  $SM$ is the unit sphere bundle and $\tau(x,\xi)$ is the length of the maximal geodesic starting at $x\in \PD M$ in the direction $\xi\in \PD_{+}SM:=\{(x,\xi)\in \PD SM: \La \nu(x),\xi\rangle\geq 0\}$.

The Fredholm-type inversion formulas of Pestov-Uhlmann are given by the following theorem:
\begin{thm}\cite[Theorem 5.4]{PU2}\label{S1:T1}
Let $(M,g)$ be a 2-dimensional Riemannian manifold. Then 
\[ f+\Wc^{2}f=\frac{1}{4\pi}\delta_{\perp}I_{1}^{*}\lb \A\sp{*}H(I_{0}f)^{-}|_{\PD_{+}SM}\rb,\quad f\in L^{2}(M).\]
\[ h+(\Wc\sp{*})^{2}h=\frac{1}{4\pi}I_{0}\sp{*}\lb\A\sp{*}H(I_{1}\Hc_{\perp}h)^{+}|_{\PD_{+}SM}\rb,\quad h\in H_{0}^{1}(M).\]
\end{thm}
Here $\Wc$ is the operator ($\Wc^{*}$ is its $L^{2}$ adjoint) on $L^{2}(M)$ defined by 
\[ \Wc f(x)=\frac{1}{2\pi}\int\limits_{S_{x}} \Hc_{\perp}\Big( \int\limits_{0}^{\tau(x,\xi)}f(\g_{x,\xi}(t))dt \Big)dS_{x}(\xi),\]
with 
\[ \Hc_{\perp}u(x,\xi)=\xi_{\perp}^{i}\lb \frac{\PD u}{\PD x^{i}}-\Gamma^{k}_{ij}\xi^{j}\frac{\PD u}{\PD \xi^{k}}\rb.\]

As shown in \cite{PU2}, for manifolds of constant curvature, $\Wc=\Wc^{*}=0$ and hence these formulas becomes exact inversion formulas.

In this note, we generalize these formulas to simple $2$-dimensional manifolds whose curvatures are close to a constant. We show that in this case, the inversion formulas are given by convergent Neumann series expansions. Our main result is a generalization of the above result:
\begin{thm}\label{S1:T2}
There exists a $C>0$ such that if $M$ is a simple 2-dimensional manifold with Gaussian curvature $K$ such that $||\nabla K||_{C^{0}}\leq C$, the following inversion formulas hold:
\[ f=(I+\Wc^{2})^{-1}\lb\frac{1}{4\pi}\delta_{\perp}I_{1}^{*}\lb \A\sp{*}H(I_{0}f)^{-}|_{\PD_{+}SM}\rb\rb,\quad f\in L^{2}(M).\]
\[ h=(I+(\Wc^{*})^{2})^{-1}\lb\frac{1}{4\pi}I_{0}\sp{*}\lb\A\sp{*}H(I_{1}\Hc_{\perp}h)^{+}|_{\PD_{+}SM}\rb\rb,\quad h\in H_{0}^{1}(M).\]
\end{thm}
Remark: The proof relies on getting bounds for the operator $\Wc$ (and hence $\Wc^{*}$) in terms of the gradient of the curvature $K$. Hence we will not give definitions of the terms appearing on the right hand side of the formulas in Theorems \ref{S1:T1} or \ref{S1:T2} which can be found in Pestov-Uhlmann's papers \cite{PU1, PU2}.

As shown in \cite{PU2}, $\Wc$ is a smoothing integral operator extendible as a map $\Wc:L^{2}(M)\to C^{\infty}(M)$ with kernel, 
\begin{equation}\label{I:K}
W(x,y)= -Q(x,\exp_{x}^{-1}(y))\frac{\left|\det(\exp_{x}^{-1})'(x,y)\right|\sqrt{g(x)}}{\sqrt{g(y)}}.
\end{equation}
The function $Q$ (equations \eqref{A:q} and \eqref{A:Q}) and the partial differential operator $\PD_{\theta}$ (equation \eqref{A:Th}) are defined in the appendix.\\

{\it Acknowledgments:} The author wishes to express his gratitude to Gunther Uhlmann and Leonid Pestov for their guidance and encouragement. 

\section{The proof}

We prove the following lemma. Here the derivatives are with respect to time. The functions $a$ and $b$ are defined in the appendix; see equation \eqref{A:SJ}.
\begin{lemma}\label{S1:L1}
Let $\vp=b\PD_{\theta}a-a\PD_{\theta}b$. Then denoting $K_{\g}=K\circ \g$, $\vp$ satisfies the following ordinary differential equation, 
\[ \vp^{(3)}+4K_{\g}\vp'+2K_{\g}'\vp=-2\PD_{\theta}K_{\g},\]
with initial conditions, 
\[ \vp(0)=\vp'(0)=\vp''(0)=0.\]
\end{lemma}
\begin{proof} 
First of all we have
\begin{equation}\label{T1:E1}
a b'-a' b\equiv 1.
\end{equation}
 For, let $\phi=ab'-ba'$. Then $\phi'=a b''-a''b$.
From equation \eqref{A:SJ} we get that $\phi'=0$ and so $\phi$ is
a constant. Since $\phi(0)=1$, we have the claim.
With this we now show that $\varphi=b\partial_{\theta}a-a\partial_{\theta}b$ satisfies the ODE above.
\[ \vp'=b'\partial_{\theta}a+b\partial_{\theta}a'-a'\partial_{\theta}b-a\partial_{\theta}b'.\]
From \eqref{T1:E1} we get 
\[b'\partial_{\theta} a+a\PD b'-b\PD_{\theta} a'-a'\PD_{\theta} b=0.\]
This gives
\[ \vp'=2(b'\PD_{\theta} a- a'\PD_{\theta} b)=2(b\PD_{\theta} a'-a\PD_{\theta}b').\]
Differentiating again, we get
\[ \vp''=2(b''\PD_{\theta} {a}+b'\PD_{\theta} a'-a''\PD_{\theta} b- a'\PD_{\theta} b').\]
Using equation \eqref{A:SJ}, this reduces to
\[ \vp''=2(-K_{\g}\varphi+b'\PD_{\theta} a'-a'\PD_{\theta}b'), \]
where $K_{\g}=K\circ\g$. Differentiating yet again, and as in the steps above, we finally get,  
\begin{equation}\label{L1:E2}
\vp^{(3)}+2{K_\gamma}'\varphi+4K_\gamma \vp'=-2\partial_{\theta}K_\gamma,\end{equation}
It now follows directly from these equations that $\vp(0)=\vp'(0)=\vp''(0)=0$.
\end{proof}

Notation: The norm $\|\cdot\|$ in the proof below denotes the $\sup$ norm unless indicated otherwise. Also in order to avoid proliferation of subscripts, we will use the same letter $C$ to denote different constants.

\begin{proof} [Proof of Theorem \ref{S1:T2}] We now prove the main theorem. We rewrite equation \eqref{L1:E2} as a first order differential equation. We get
\begin{equation*}
\left(
\begin{matrix}
\vp_{1}'\\
\vp_{2}'\\
\vp_{3}'
\end{matrix}
\right)=
\left(
\begin{matrix}
0 & 1 & 0\\
0 & 0 & 1\\
-2K_{\g}' & -4K_{\g} & 0
\end{matrix}
\right)
\left(
\begin{matrix}
\vp_{1}\\
\vp_{2}\\
\vp_{3}
\end{matrix}
\right)+
\left(
\begin{matrix}
0\\
0\\
-2\partial_{\theta}K_{\g}
\end{matrix}
\right)
\end{equation*}
where
\[\varphi_{1}=\varphi, \varphi_{2}=\varphi' \text { and } \varphi_{3}=\varphi''.\]
For simplicity, let us write this as a system of the form
\[ X'(t)=A(t) X(t)+B(t),\]
where $X,B$ and matrix $A$ depend also on
$(x,\xi)$). From \cite{CL}, since $X(0)=0$, we have a solution of this
differential equation to be
\[ X(t)= \Phi(t)\int\limits_{0}^{t}\Phi^{-1}(s)B(s)\,ds.\]
where $\Phi$ is the fundamental matrix of the homogeneous differential equation,
\[X'(t)=A(t)X(t).\]
Since the manifold is compact, we have $\|\Phi\|, \|\Phi^{-1}\|<\infty$.
From the relation,
\[ \PD_{\theta} K_\gamma=( \xi_{\perp},\nabla K_\gamma) , \]
and using the fact that $S M$ is compact, we have a $C$ such that
\[ |\PD_{\theta} K_\gamma|\leq C \|\nabla K\|.\]
Combining these inequalities we get,
\[ |\vp'(x,\xi,t)|\leq |X(x,\xi,t)|\leq C t \|\nabla K\|.\]
for some $C>0$. Since
\[ \varphi(t)=\int\limits_{0}^{t}\vp'(s)ds,\] we have
\[ |\varphi(t)|\leq Ct^{2}\|\nabla K\|.\] 
We can initially work with $\vp''(t)$ and by the same  arguments as above, we also get 
\[ |\varphi(t)|\leq Ct^{3}\|\nabla K\|,\]
for a different constant $C$.

Since the manifold is simple, we have $b\neq 0$ for $t\neq 0$, since $b(0)=0$. Now we write $b(t,x,y)=t\tilde{b(x,y,t)}$ with
$\tilde{b}\neq 0.$ Therefore for a suitable $C>0$,
\[ |q(x,\xi,t)|=|\partial_\theta \frac{a}{b}|\leq Ct\|\nabla K\|,\]
where the norm of $\nabla K$ is the $\sup$ norm. Since $tQ(x,t\xi)=q(x,\xi,t)$, we have 
\[ |Q(x,t\xi)|\leq C\|\nabla K\|.\]
Since the remaining terms in 
\[ W(x,y)= -Q(x,\exp_{x}^{-1}(y))\frac{\det(\exp_{x}^{-1})'(x,y)\sqrt{g(x)}}{\sqrt{g(y)}}\]
are bounded above by compactness of $M$, we have 
\[ \|W\|\leq C\|\nabla K\|.\]
Therefore we have 
\[ \|\Wc\|_{L^{2}\to C^{\infty}(M)}\leq C\|\nabla K\|.\]
So now choosing $\|\nabla K\|$ to be small enough, we have $\|\Wc\|<1$. Hence we have inversion formulas involving Neumann series expansions recovering the function from its geodesic ray transform. A similar argument works for the recovery of the solenoidal part of a vector field from its geodesic ray transform. This completes the proof of the theorem.
\end{proof}
\appendix
\section{The kernel of $\Wc$}
For completeness and because the function $q$ defined in equation \eqref{A:q} is critical for the proof of Theorem \ref{S1:T2}, we sketch below, Pestov-Uhlmann's \cite{PU2} derivation of the integral kernel of the operator $\Wc$.

Recall that the operator $\Wc$ is defined as  
\[ \Wc f(x)=\frac{1}{2\pi}\int\limits_{S_{x}M}\Hc_{\perp}\int\limits_{0}^{\tau(x,\xi)}f(\g(x,\xi,t))dt\,dS_{x},\]
where 
\[ \Hc_{\perp}=\xi_{\perp}^{i}\lb \frac{\PD}{\PD x^{i}}-\Gamma_{ij}^{k}\xi^{j}\frac{\PD}{\PD \xi^{k}}\rb.\]
For a function $u$ on $SM$, $\frac{\PD}{\PD \xi^{k}}u$ is defined as 
\[ \frac{\PD}{\PD \xi^{k}}u=\frac{\PD}{\PD \xi^{k}}(u\circ p)|_{|\xi|=1}, \text{ where } p(x,\xi)=(x,\xi/|\xi|).\]
We have 
\[ \Wc f(x)=\frac{1}{2\pi}\int\limits_{S_{x}M}\int\limits_{0}^{\tau(x,\xi)}\La\nabla f,\Hc_{\perp}\g\rangle.\]
We define two Jacobi vector fields along the geodesic $\g(x,\xi,t)$ as follows: Let $x(s), -\ve<s<\ve$ be a curve starting at $x$ in the direction $\xi_{\perp}$. Now parallel translate the vector $\xi$ along this curve, call it $\xi(s)$ and consider the variation by geodesics, $\g(x(s),\xi(s),t)$. The vector field
\[X(x,\xi,t)=\frac{d}{ds}|_{s=0}\g(x(s),\xi(s),t),\]
is a Jacobi vector field along $\g$ with the following initial conditions, 
\[ X(x,\xi,0)=\xi_{\perp},\quad D_{t}X(x,\xi,0)=0.\]
It can be also be written as 
\[ X(x,\xi,t)=\Hc_{\perp}\g(x,\xi,t).\]
We now define another Jacobi vector field by considering the variation by geodesics, $\g(x,\xi(s),t)$, where $\xi(s)$ is a smooth curve in $S_{x}M$ with initial tangent vector $\xi_{\perp}$. The Jacobi vector field
\begin{equation}\label{A:Th}
\PD_{\theta}\g(x,\xi,t):= Y(x,\xi,t)=\frac{d}{ds}|_{s=0} \g(x,\xi(s),t)
\end{equation}
has initial conditions,
\[ Y(x,\xi,0)=0,\quad D_{t}Y(x,\xi,0)=\xi_{\perp}.\]
Since $X$ and $Y$ are vector fields normal to $\g$ and because of dimensional reasons these two fields must be proportional to the parallel translate of the vector $\xi_{\perp}$ along the geodesic $\g$. Let this parallel translate be denoted $\dot{\g}_{\perp}$. Then there exists two smooth functions $a(x,\xi,t)$ and $b(x,\xi,t)$ such that 
\[ X=a\dot{\g}_{\perp},\quad Y=b\dot{\g}_{\perp}.\]
The functions $a$ and $b$ satisfy the scalar Jacobi equations,
\begin{equation}\label{A:SJ}
a''+Ka=b''+Kb=0,
\end{equation}
with 
\[ a(x,\xi,0)=1,\quad a'(x,\xi,0)=0,\quad b(x,\xi,0)=0,\quad b'(x,\xi,0)=1.\]
We now write, 
\begin{align*}
 \Wc f(x)&=\frac{1}{2\pi}\int\limits_{S_{x}M}\int\limits_{0}^{\tau(x,\xi)}\La\nabla f,\Hc_{\perp}\g\rangle dt dS_{x}\\
 &=\frac{1}{2\pi}\int\limits_{S_{x}M}\int\limits_{0}^{\tau(x,\xi)}\frac{a}{b}\La\nabla f,Y\rangle dt dS_{x}\\
  &=-\frac{1}{2\pi}\int\limits_{0}^{\tau(x,\xi)}\int\limits_{S_{x}M}\PD_{\theta}(\frac{a}{b})f\circ \g dS_{x}dt.
 \end{align*}
 We now define a function $q$ on 
 \[ G=\{(x,\xi,t):(x,\xi)\in SM, -\tau(x,-\xi)<t<\tau(x,\xi), t\neq 0\}\]by 
 \begin{equation}\label{A:q}
q(x,\xi,t)=\PD_{\theta}\lb\frac{a}{b}\rb.
\end{equation}
Using this we define a function $Q\in C^{\infty}(TM)$ by, 
\begin{equation} \label{A:Q}
Q(x,t\xi)=tq(x,\xi,t).
\end{equation}
The existence of this function $Q$ follows from the fact that the geodesic $\g(x,\xi,t)$ is smooth as a function of $(x,t\xi)$ and the initial conditions of the Lemma \ref{S1:L1}. We now get the integral kernel $W(x,y)$ in equation \eqref{I:K} by a change of variables involving the  inverse of the exponential map. 
\bibliographystyle{alpha}
\bibliography{mybibfile}
\end{document}